\documentclass[conference,dvipdfmx]{IEEEtran}
\IEEEoverridecommandlockouts
\usepackage{amsmath,amssymb,amsfonts,amsthm}
\usepackage{graphicx}

\usepackage{mathrsfs}
\usepackage{cases}
\usepackage{citesort}
\usepackage{balance}

\newcommand{\argmax}{\mathop{\rm arg~max}\limits}
\newcommand{\argmin}{\mathop{\rm arg~min}\limits}
\newcommand{\minimize}{\mathop{\rm minimize\ \ }\limits}

\newcommand{\prox}{\mathrm{prox}}

\newcommand{\fix}{{\mathrm {Fix}}}

\renewcommand{\baselinestretch}{0.9903}

\newtheorem{example}{Example}
\newtheorem{remark}{Remark}
\newtheorem{definition}{Definition}
\newtheorem{theorem}{Theorem}

\def\BibTeX{{\rm B\kern-.05em{\sc i\kern-.025em b}\kern-.08em
    T\kern-.1667em\lower.7ex\hbox{E}\kern-.125emX}}
\begin{document}

\title{
  A Hierarchical Convex Optimization for  Multiclass SVM Achieving Maximum Pairwise Margins \\ with Least Empirical Hinge-Loss
}

\author{\IEEEauthorblockN{Yunosuke Nakayama,  Masao Yamagishi, and Isao Yamada}
  \IEEEauthorblockA{\textit{Department of Information and Communications Engineering},
    \textit{Tokyo Institute of Technology},
Japan \\
Email: \{nakayama, myamagi, isao\}@sp.ce.titech.ac.jp}
}

\maketitle

\begin{abstract}
  In this paper, we formulate newly a hierarchical convex optimization for  multiclass SVM achieving maximum 
pairwise margins with least empirical hinge-loss. 
This optimization problem is a most faithful as well as robust multiclass extension of an NP-hard hierarchical optimization appeared for the first time in the seminal paper by  C.~Cortes and V.~Vapnik almost 25 years ago.
 By extending the very recent fixed point theoretic idea [Yamada-Yamagishi 2019] with 
the generalized hinge loss function [Crammer-Singer 2001], we show that {\sl the hybrid steepest descent method} [Yamada 2001] in the computational fixed point theory is applicable to this much more complex hierarchical convex optimization problem. 
\end{abstract}


\section{Introduction}
\label{sec:intro}
For the classical two-group classification problem,  {\sl the  soft-margin hyperplane} (or {\sl the soft-margin SVM}) was introduced in the seminal paper \cite[Sec.3]{cortes1995support} by C.~Cortes and V.~Vapnik 
as the solution of a naive convex relaxation of a certain NP-hard optimization problem  (called, in this paper,  Cortes-Vapnik problem)\footnote{
                            In this paper, for simplicity, we focus on the linear SVM because the nonlinear SVM 
                            (and its extensions for multiclass classification) 
                            exploiting the so-called Kernel trick can
                            be viewed as an instance of the linear classifiers in the \textit{Reproducing Kernel Hilbert Spaces} (RKHS).
                            }.  The Cortes-Vapnik problem has a  hierarchical structure, for finding a special hyperplane (or linear classifier) with maximum margin among all hyperplanes achieving minimum number of misclassified training samples\footnote{Since minimizing the number of misclassified training samples is known to be NP-hard, \cite[Sec.3]{cortes1995support} proposes to minimize instead 
$\phi(\mathbf{w},b)$ defined as {\sl the sum of deviations}  of training errors.   
This idea was certainly a step ahead of the nowadays standard computational techniques known as {\sl the $\ell_{1}$ based convex relaxations} in the sparsity aware data sciences.}, where we observe not only (i) the number of misclassified samples as the 1st stage optimization criterion, but also (ii)  the margin as the 2nd stage optimization criterion for the special hyperplane. The solution of Cortes-Vapnik problem is clearly an ideal extension of the {\sl optimal separating hyperplane} \cite{vapnik82,vapnik1998statistical} (one of the most monumental landmarks in the history of pattern recognition),  well-defined only for linearly separable training data but achieving maximum margin without misclassified samples, to the general case where the training data is possibly linearly nonseparable. 

For the advancement of pattern recognition and its applications,  
it is important  to establish a reliable computational algorithm for approximating  as faithfully as possible the solution of the Cortes-Vapnik problem.   
We believe that such a most faithful approximation of  the solution of the Cortes-Vapnik problem can be achieved by solving {\sl the hierarchical convex relaxation} \cite[16.4.2] {yamada2017hierarchical} which is formulated with 
replacement of only the 1st stage optimization criterion in the Cortes-Vapnik problem by  {\sl the sum of deviations} of training errors while keeping its hierarchical structure.   However the reliable numerical solution of  the hierarchical convex relaxation of 
the Cortes-Vapnik problem had  not been available until \cite{yamada2017hierarchical}    
 since \cite[Sec.3]{cortes1995support} employed 
 a naive convex relaxation\footnote{
                                    The function $\phi (\mathbf{w},b)$ can be expressed with a \textit{hinge loss function}.  
                                    The squared margin of the linear classifier with $(\mathbf{w},b)$ is given by $ 2^{-1}\psi^{-1}(\mathbf{w},b) = \frac{1}{\|\mathbf{w}\|^2}$.
                                    }:
    \begin{align}
        \minimize_{(\mathbf{w},b)\in \mathbb{R}^N \times \mathbb{R}} \psi(\mathbf{w},b) + C \phi(\mathbf{w},b)
        \label{pro:naive},
    \end{align}
        with a \textit{tuning parameter} $C > 0$,  just for minimizing the weighted sum of 
        $\psi(\mathbf{w},b)$ and $\phi(\mathbf{w},b)$, and put  the computational difficulty, caused by the hierarchical structure, in the hand of the delicate tuning parameter $C$
 \footnote{Note that the suggestion, in \cite[Sec.3]{cortes1995support}, of using a sufficiently large $C$ does not provide us with any practical computational strategy in general case where  the target is linearly nonseparable data. This is because in principle we can not judge whether the minimization of {\sl the sum of deviations} of training errors is achieved or not before convergence of any iterative solver for convex optimization with finite weight. 
 See also  Remark~\ref{rem:soft-marginSVM}(a) for  the awkwardness in the formulation (\ref{pro:naive}) of the soft-margin SVM.}.
Very recently, by applying {\sl the hybrid steepest descent method}  \cite{hybrid} in the computational fixed point theory 
to special nonexpansive operators designed through the art of {\sl proximal splitting} in convex analysis, \cite{yamada2017hierarchical} succeeded in establishing a reliable numerical algorithm to solve {\sl the hierarchical convex relaxation} of the Cortes-Vapnik problem for the binary classification. 
  
Various extensions of the native convex relaxation, in (\ref{pro:naive}) of the Cortes-Vapnik problem,  
to be applicable to multiclass classification,  
    have also been reported by many researchers   
    (see, e.g.,\cite{blanz1996comparison,kressel1999pairwise,crammer2001algorithmic,weston1998multi,guermeur2002combining,tatsumi2010multiobjective}).
Such extensions include \textit{combination approaches} and \textit{direct approaches}.
The combination approach applies the binary soft-margin SVM to multiple independent binary classification problems 
    defined for each pair of disjoint sets of classes in the original multiclass problem \cite{blanz1996comparison,kressel1999pairwise}.
Since the combination approaches cannot capture the correlations between the different classes \cite{crammer2001algorithmic},   the direct approach has been formulated, as 
 a single convex optimization problem,   by extending in various ways \cite{vapnik1998statistical,crammer2001algorithmic,weston1998multi,guermeur2002combining,tatsumi2010multiobjective} 
 the cost functions  $\psi(\mathbf{w},b)$ and $\phi(\mathbf{w},b)$ in  (\ref{pro:naive}). 
 However, as  pointed out in \cite{tatsumi2010multiobjective},   the existing extensions especially for the margin function 
        have not yet succeeded in robustifying fairly for every pair of classes in the multiclass problem.
Moreover, the hierarchical convex relaxation \cite{yamada2017hierarchical} of the Cortes-Vapnik problem 
has never been extended to  the multiclass problem.  

In this paper, for an ideal multiclass extension of SVM,   
we newly formulate  a hierarchical convex optimization problem 
for finding a special multiclass linear classifier with maximum pairwise margins among all multiclass linear classifiers 
achieving minimum sum of deviations for linearly nonseparable data 
where the criterion in the 2nd stage optimization is designed to robustify fairly for every pair of classes through 
maximization of the least one among all 
pairwise margins.  Moreover, by extending the fixed point theoretic characterization in 
\cite[16.4.2]{yamada2017hierarchical} with 
the generalized hinge loss function \cite{crammer2001algorithmic}, we show that the hybrid steepest descent method 
is applicable to this much more complex hierarchical convex optimization problem. 
Numerical experiments demonstrate that the proposed classifier achieves larger smallest pairwise margin, with least empirical hinge loss, than the existing multiclass extensions of binary SVM.  

\section{PRELIMINARIES}
%
%
\subsection{Multiclass Support Vector Machine}
We consider a multiclass classification
    with a given training dataset: 
        \begin{align}
            \mathcal{D} := & \ \{(\mathbf{x}_i, y_i) \in \mathbb{R}^N \times \mathcal{I}_K \mid i=1,2,\dots,M\} \label{eq:training_data} \\
                \mathcal{D}_j := & \ \{ \mathbf{x}_i \in \mathbb{R}^N \mid (\mathbf{x}_i,j) \in \mathcal{D}  \} \quad(j \in \mathcal{I}_K),
        \end{align} 
    where 
    \( y_i \in \mathcal{I}_K := \{1,2,\dots,K\} \) stands for the label assigned to
        \(\mathbf{x}_i\).
The multiclass linear classifier is a mapping:
        \begin{align}
            \mathscr{L}_{\mathbf{p}} :\mathbb{R}^N \to \mathcal{I}_K : 
            \mathbf{x} \mapsto 
            \min \left(\argmax_{j \in \mathcal{I}_K} \mathbf{w}_j^\top \mathbf{x} + b_j\right)
        \end{align}
    defined with \(\mathbf{p} := (\mathbf{w}_j,b_j)_{j\in\mathcal{I}_K} \in \mathbb{R}^{(N+1)K} \).
In this paper, the dataset \(\mathcal{D}\) is said to be \textit{linearly separable by a multiclass machine}
    if there exists \(\mathbf{p} \in \mathbb{R}^{(N+1)K}\) satisfying
        \begin{align}
        (\forall i \in \{1,2,\dots,M\}) \ 
        \argmax_{j \in \mathcal{I}_K} \mathbf{w}_j^\top \mathbf{x}_i + b_j = \{y_i\}.
        \label{condition1}
        \end{align} 
To design suitable \(\mathbf{p}\) for a not necessarily separable dataset \(\mathcal{D}\),
    many extensions \cite{crammer2001algorithmic,chierchia2015proximal,tatsumi2010multiobjective} 
     of the binary support vector machine (SVM) \cite{vapnik1998statistical,cortes1995support}
        have been proposed for the task of multiclass classification.

In the following, along a geometrical point of view found in \cite{yamada2017hierarchical},
    we introduce a criterion, used in \cite{crammer2001algorithmic,chierchia2015proximal} to optimize \(\mathbf{p}\).
This criterion can be interpreted as a convex relaxation of the number of misclassified samples 
    by the linear classifier \(\mathscr{L}_{\mathbf{p}}\)with tunable margins. 

To extend naturally the binary SVM for multiclass classifier \(\mathscr{L}_{\mathbf{p}}\), we follow the strategy of the binary SVM.
By the definition of argmax,
    the classifier \(\mathscr{L}_{\mathbf{p}}\) can be built from the binary classifiers for every pair  
        \((r,s)\in \mathcal{I}_K\times \mathcal{I}_K,\textrm{\ s.t.}\ r\neq s\):
            \begin{align}
                \hspace{-2mm}
                \mathcal{L}_{(\mathbf{p},r,s)} : \mathbb{R}^N \to \{r,s\} : 
                \mathbf{x} \mapsto 
                \begin{cases}
                    r \quad \textrm{if} \ \pmb{\omega}_{rs}^\top \mathbf{x} + \beta_{rs}  > 0,  \vspace{2mm}\\
                    s \quad \textrm{if} \  \pmb{\omega}_{rs}^\top \mathbf{x} + \beta_{rs} \leq 0,  
                \end{cases}
                \hspace{-2mm}\label{def:6}
            \end{align}
    where \(\pmb{\omega}_{rs} := \mathbf{w}_r - \mathbf{w}_s\) and \(\beta_{rs} := b_r - b_s\)
        (Note: \( \mathcal{L}_{(\mathbf{p},r,s)} = \mathcal{L}_{(\mathbf{p},s,r)}\) holds for all \(r\neq s\)).
Note that if \(\mathcal{D}\) is linearly separable,
    there also exists infinitely many \(\mathbf{p} \in \mathbb{R}^{(N+1)K} \) satisfying
        for every pair of \(r,s \in \mathcal{I}_K\times \mathcal{I}_K,\textrm{\ s.t.}\ r\neq s\)
            \begin{align}
                \hspace{-3mm}
                \left.
                \begin{array}{ll}
                \mathcal{D}_r \subset \Pi_{(\mathbf{p},r,s)}^{[1,\infty)} := \{\mathbf{x}\in\mathbb{R}^N \mid \pmb{\omega}_{rs}^\top\mathbf{x} + \beta_{rs} \geq 1 \}\\
                \mathcal{D}_s \subset \Pi_{(\mathbf{p},r,s)}^{(-\infty,-1]} := \{\mathbf{x}\in\mathbb{R}^N \mid \pmb{\omega}_{rs}^\top\mathbf{x} + \beta_{rs} \leq -1 \}.
                \end{array}
                \right \}\label{con:7}
            \end{align}
The closed half-spaces 
        \(\Pi_{(\mathbf{p},r,s)}^{[1,\infty)}\) and \(\Pi_{(\mathbf{p},r,s)}^{(-\infty,-1]}\)
    are main players in the following consideration on \(\mathscr{L}_{\mathbf{p}}\)
    even for not linearly separable data \(\mathcal{D}\).
In this paper,
    the \textit{pairwise margin} of \(\mathcal{L}_{(\mathbf{p},r,s)}\) in (\ref{def:6}) is defined by 
        \begin{align}
            \textrm{dist} \left(\Pi_{(\mathbf{p},r,s)}^{[1,\infty)} ,\Pi_{(\mathbf{p},r,s)}^{0}  \right)
            =
            \frac{1}{\|\pmb{\omega}_{rs}\|} \left( = \frac{1}{\|\mathbf{w}_r - \mathbf{w}_s\|} \right),
            \label{margin}
        \end{align}
    where 
        \(\|\cdot \| \) stands for the Euclidean norm in \(\mathbb{R}^N\)
  and
        \(\textrm{dist} \left(\Pi_{(\mathbf{p},r,s)}^{[1,\infty)} ,\Pi_{(\mathbf{p},r,s)}^{0} \right) \)
    stands for the distance between the closed half-space and the decision hyperplane. 
By using the distance between 
        \(\mathbf{x} \in \mathbb{R}^N \)
    and
        \(\Pi_{(\mathbf{p},r,s)}^{[1,\infty)}\)
    given by
        \begin{align}
            {\rm d}\left({\bf x},\Pi_{(\mathbf{p},r,s)}^{[1,\infty)}\right) =
            \begin{cases}
            \frac{1- (\pmb{\omega}_{rs}^{\top} {\bf x} + \beta_{rs})}{\|\pmb{\omega}_{rs}\|}
                &{\rm if} \ {\bf x} \notin \Pi_{(\mathbf{p},r,s)}^{[1,\infty) },\\
            0 &{\rm otherwise},
            \end{cases}
        \end{align}
    we deduce 
      \begin{align}
        \hspace{-3mm}
        \Phi_{\mathcal{D}} (\mathbf{p}) 
        &:=  \sum_{i=1}^M  \tilde\delta_{\mathbf{p}} (i)\\
        &=
        \sum_{i=1}^M  \max \left\{0,\max_{s\in\mathcal{I}_K\setminus\{y_i \}} 
         [1- (\pmb{\omega}_{y_is}^\top\mathbf{x}_i + \beta_{y_is})] \right\},
          \hspace{-1.5mm}
          \label{hinge loss}
       \end{align}
    where
        \begin{align}
            \tilde\delta_{\mathbf{p}} (i) := 
            \max_{s\in\mathcal{I}_K\setminus\{y_i \}}
            \|\pmb{\omega}_{y_is} \|
            {\rm d}\left(\mathbf{x}_i,\Pi_{(\mathbf{p},y_i,s)}^{[1,\infty)}\right).
        \end{align}
The right hand side of (\ref{hinge loss}) called as 
    \textit{(generalized) hinge loss function} \cite{chierchia2015proximal,wang20071,liu2006multicategory} 
        introduced originally in \cite{crammer2001algorithmic}.
For each \(i \in \{1,2,\dots,M\}\),
     we see that
     \begin{align}
        \tilde\delta_{\mathbf{p}}(i) > 0
        \Leftrightarrow 
        \ (\exists s \in \mathcal{I}_K\setminus\{y_i\})\ 
    \mathbf{x}_i \notin \Pi_{(\mathbf{p},y_i,s)}^{[1,\infty)}
    \label{hinge loss2}
     \end{align}
    and therefore the hinge loss function in (\ref{hinge loss}) can be seen as a convex relaxation of 
        \begin{align}
            \mathscr{E}_{\mathbf{p}}(\mathcal{D})
            \hspace{-1mm}:= \hspace{-1mm}
            \left |  \left \{i \in \{1,\dots,M\} \left| 
            \exists s \in \mathcal{I}_K
            \hspace{-1mm}\setminus \hspace{-1mm}\{y_i\},  
            \mathbf{x}_i \in \Pi_{(\mathbf{p},y_i,s)}^{(-\infty,1)} \right.
        \right \} \right | \notag
            \end{align}
    which is desired to be minimized because it is 
     the  total number of misclassified samples \(\mathbf{x}_i\)
     together with samples \(\mathbf{x}_i\) likely to be misclassified, by \(\mathscr{L}_{\mathbf{p}}\),
        satisfying
        \begin{align}
            \mathbf{x}_i \in  \Pi_{(\mathbf{p},y_i,s)}^{(-\infty,1)} := \{\mathbf{x}\in\mathbb{R}^N \mid   \pmb{\omega}_{y_is}^\top\mathbf{x} + \beta_{y_is} < 1 \}.
        \end{align}
Among all \(\mathbf{p} \in \mathbb{R}^{(N+1)K}\) of the same \(\mathscr{E}_{\mathbf{p}}(\mathcal{D})\),
    in order to achieve higher generalization performance,
        it is preferred to choose \(\mathbf{p}\) achieving 
            larger pairwise margin \(\textrm{dist} \left(\Pi_{(\mathbf{p},r,s)}^{[1,\infty)} ,\Pi_{(\mathbf{p},r,s)}^{0} \right) \)
                for every pair \((r,s) \in \mathcal{I}_K \times \mathcal{I}_K, \textrm{\ s.t.\ }, r\neq s\)
                    because, for unknown \(\mathbf{x}\) satisfying \(\mathscr{L}_{\mathbf{p}}(\mathbf{x}) = r\),
                        the risk for misclassification \(\mathscr{L}_{\mathbf{p}}(\mathbf{x+\pmb \epsilon}) = s\)
                            is desired to be suppressed even after being contaminated by noise \(\pmb \epsilon \).  
 
For such a purpose, \cite{crammer2001algorithmic} proposed for general training dataset \(\mathcal{D}\) in \eqref{eq:training_data}
 the following design of multiclass linear classifier:
    \begin{align}
                \minimize_{\mathbf{p}\in\mathbb{R}^{(N+1)K}}
            \frac{1}{2}  \sum_{r,s\in \mathcal{I}_K \ (r<s)} \| \mathbf{w}_r - \mathbf{w}_s\|^2 
                + C  \Phi_{\mathcal{D}} (\mathbf{p}) ,
        \label{pro:Crammer_Singer}
    \end{align}
        with a tuning parameter \(C > 0\)
        under the simplified assumption \(b_j = 0 \  ( j \in \mathcal{I}_K) \).        
        \begin{remark}
          \label{rem:soft-marginSVM}
\begin{enumerate}
    \item[(a)]
\textit{The solution of (\ref{pro:Crammer_Singer}) is clearly a direct extension of the binary soft-margin SVM \cite{cortes1995support}
which has been utilized widely as a standard binary-linear classifier.
As seen in \cite{cortes1995support},
    the binary soft-margin SVM is the solution of a naive convex relaxation of a certain NP-hard hierarchical optimization (Cortes-Vapnik problem) \cite{cortes1995support}
    by the complete loss of any hierarchical structure (see also \cite{yamada2017hierarchical}).
However, this naive convex relaxation has no guarantee to reproduce
 the original SVM \cite{vapnik1998statistical} 
    established specially for linearly separable training dataset.
To overcome this awkwardness, 
 a novel hierarchical convex relaxation has been established in \cite{yamada2017hierarchical}.
}
\item[(b)]
\textit{Certainly, for the robustness against noise,
    the risk of misclassification by the multiclass linear classifier is desired to be suppressed 
    for every pair of different classes.
However in the existing extensions \cite{crammer2001algorithmic,chierchia2015proximal} of 
        the binary soft-margin SVM \cite{cortes1995support},
        only the average of the squared inverse of all pairwise margins has been suppressed (see (\ref{pro:Crammer_Singer}))
        but any special care for the most risky pair 
        has not been taken strategically.
        }
      \end{enumerate}
    \end{remark}
Motivated by these facts,  
 we will present a practical hierarchical formulation to extend the SVM for multiclass classification in Section~\ref{sec:pagestyle}.
In the following,
  we will give key mathematical tools for the hierarchical convex optimization.
%
%
\subsection{Hierarchical Convex Optimization with Proximal Splitting Operator} 
Let $\mathcal{H}$ be a real Hilbert space.
For a nonexpansive operator $T : \mathcal{H} \to \mathcal{H}$,
    i.e.,
    $ \| T(x) - T(y) \|_{\mathcal{H}} \leq \| x - y\|_{\mathcal{H}}\ (\forall x,y \in \mathcal{H}),$ 
        the set of all fixed points of  $T$, denoted by $\fix(T) := \{ x \in \mathcal{H} \mid T(x) = x\}$,
            is known to be a closed convex set.
If $T_i : \mathcal{H} \to \mathcal{H}\ (i = 1,2,\dots,m)$ are averaged nonexpansive 
    with $\cap_{i=1}^m \fix(T_i) \neq \varnothing$, 
        $T_m \circ T_{m-1}\circ \cdots \circ T_1$ is also averaged nonexpansive and $\fix(T_m \circ T_{m-1}\circ \cdots \circ T_1) = \cap_{i=1}^m \fix(T_i)$
            (See \cite[Sec.4.5]{bauschke2011convex} for averaged nonexpansive operators).
Moreover,
    if a nonexpansive operator \( T \) with \(\fix(T) \neq \varnothing\) and
        the gradient \(\nabla \Theta : \mathcal{H} \to \mathcal{H}\) of an L-smooth convex function \(\Theta : \mathcal{H} \to \mathbb{R}\) are computable,
 we can minimize \(\Theta\)  
    over \(\fix(T) \) by the \textit{Hybrid Steepest Descent Method (HSDM)} \cite{yamada2017hierarchical,hybrid}:
        \begin{align}
            \mathbf{z}_{n+1} = T(\mathbf{z}_n) - \lambda_{n+1} \nabla \Theta (T(\mathbf{z}_n))
            \label{HSDM}
        \end{align}
    with a slowly vanishing sequence $(\lambda_n)_{n\geq1}$ $\subset [0,\infty)$,
    under reasonable conditions 
    (see, e.g., \cite{yamada2017hierarchical,hybrid,ogura2003nonstrictly,yamada2011minimizing} for the technical detail of the method).
     
\textit{The proximal splitting operators},
    developed in the art of \textit{proximal splitting techniques},
    are nonexpansive operators designed with  the so-called 
    \textit{proximity operator}\footnote{For \(\varphi \in \Gamma_0(\mathcal{X})\),
        i.e.,
        \(\varphi\) is a proper lower semicontinuous convex function defined on a real Hilbert space \(\mathcal{X}\),
        the proximity operator of \(\varphi\) is defined as
            \( \prox_{\varphi}: \mathcal{X} \to \mathcal{X} : x \mapsto \argmin_{y\in\mathcal{X}}
                 \left[\varphi(y) + \frac{1}{2} \| y-x\|_{\mathcal{X}}^2\right]\).
For a closed convex set $C \subset \mathcal{X}$,
    the proximity operator of 
    \( \iota_C : \mathcal{X} \to (-\infty,\infty]
    : x \mapsto 
    \begin{cases}
        0 &\textrm{if\ }x \in C, \\
        \infty &\textrm{otherwise},
        \end{cases}
        \)
        is 
    given by the metric projection:
\(
    P_C : \mathcal{X}\to\mathcal{X} : x \mapsto \argmin_{y\in C}\|x-y\|_{\mathcal{X}}.
\)
If $\prox_{\varphi}$ is available  as a computable operator,
    $\varphi \in \Gamma_0(\mathcal{X})$ is said to be \textit{proximable}.
    }
    as their computable building blocks.
Proximal splitting operators
    are useful to characterize \(\mathcal{S}_\star := \argmin (f+g)(\mathcal{X}) \neq \varnothing \)
        for proximable \(f,g \in \Gamma_{0}(\mathcal{X})\) in terms of their fixed-point sets.
        \begin{example}[DRS operator \cite{yamada2017hierarchical,eckstein1992douglas})]
Let $f,g \in \Gamma_0(\mathcal{X})$ and \(\argmin(f+g)(\mathcal{X}) \neq \varnothing \).
Then    
    \textit{the DRS operator} 
        \begin{align}
        T_{\textrm{DRS}} := (2\prox_{f} - \mathrm{Id}) \circ (2\prox_{g} - \mathrm{Id}),
        \label{TDRS}
        \end{align}
    where $\mathrm{Id}$ stands for the identity operator, is nonexpansive and satisfies
        \begin{align}
            \prox_{g}(\fix(T_{\textrm{DRS}})) = \argmin(f+g)(\mathcal{X}).
            \label{eq:TDRS}
        \end{align}
      \end{example}
      As introduced in \cite{yamada2017hierarchical,yamada2011minimizing},
    plugging the proximal splitting operators into the hybrid steepest descent method
    has a great deal of potential in various applications of the hierarchical convex optimization: 
    \begin{align}
      \left.
      \begin{array}{ll} 
        \minimize \Psi(x^\star)\\ 
        {\rm subject\ to} \hspace{2mm}
                           x \in \mathcal{S}_\star := \argmin_{x \in \mathcal{X}} \Phi(x) \neq \varnothing,
      \end{array}
     \right \} \label{HO}
    \end{align}
    where \(\Psi, \Phi \in \Gamma_{0}(\mathcal{X})\).
\section{Robust Hierarchical Convex Multiclass  SVM}
\label{sec:pagestyle}
\subsection{Problem Formulation for rHC-mSVM}
    We propose Robust Hierarchical Convex multiclass SVM (rHC-mSVM) for multiclass classification
    as the solution of the following hierarchical convex optimization problem.
    \begin{definition}[\textit{Robust Hierarchical Convex multiclass SVM}]
      \label{def:rHC-mSVM}
For a given training dataset \(\mathcal{D}\) in (2)
    and \( \Phi_{\mathcal{D}}\) in (10)-(12),
    assume \(\mathfrak{S}_\star := \argmin \Phi_{\mathcal{D}}  (\mathbb{R}^{(N+1)K}) \neq \varnothing \).
Then, the rHC-mSVM is defined as a solution,
    say \(\mathbf{p}^{\star\star} \in \mathfrak{S}_\star\),
    of
    \begin{align}
        \hspace{-4mm}
        \left.
        \begin{array}{ll}
        \minimize 
        {\displaystyle \max_{r,s\in \mathcal{I}_K \ (r<s)} \| {\bf w}^\star_r -  {\bf w}^\star_s\|}
        (=: \Psi(\mathbf{p^\star}) )
        \\ {\rm subject\ to\ }
        \mathbf{p}^\star :=  (\mathbf{w}_j^\star,b_j^\star)_{j\in\mathcal{I}_K}  \in 
        \mathfrak{S}_\star.
        \end{array}
        \right \}
        \label{Proposition1}
    \end{align}
  \end{definition}
  \begin{remark}
    \label{rem:rHC-mSVM}
    \begin{enumerate}
    \item[(a)]
    \textit{Unlike the existing formulation (\ref{pro:Crammer_Singer}),
                Problem (\ref{Proposition1}) has a hierarchical structure
                where the criterion for the first stage optimization
                is (generalized) hinge loss function \(\Phi_{\mathcal{D}}\).
            As seen from (\ref{hinge loss2}),
                  if  \(\mathcal{D}\) is linearly separable,
                    every \(\mathbf{p}^\star \in \mathfrak{S}_\star\) achieves 
                    \(\Phi_{\mathcal{D}}(\mathbf{p}^\star) = 0\)
                and 
                    \( \mathscr{E}_{\mathbf{p}^\star}(\mathcal{D}) = 0 \).
            In particular for \(K=2\),
              since
                \(\Psi(\mathbf{p}^\star)\) 
                 becomes the inverse of the margin of the binary classifier
                introduced in \cite{vapnik1998statistical},
                the solution of (\ref{Proposition1}) can be seen
                as a multiclass extension of \cite[Sec. 16.4.2]{yamada2017hierarchical}
                and therefore a natural extension of \cite{vapnik1998statistical}
                to general dataset \(\mathcal{D}\)
                not necessarily linearly separable.
            On the other hand,
                even if \(K=2 \) and \(\mathcal{D}\) is linearly separable, 
                the solution of (\ref{pro:Crammer_Singer})
                cannot reproduce, in general,
                the binary classifier in \cite{vapnik1998statistical}.
         } 
    \item[(b)]
    \textit{To the best of the authors' knowledge,
                Problem (20) presents 
                the first design strategy
                of the multiclass linear classifier 
                which 
                maximizes
                all pairwise margins uniformly,
                i.e.,
                maximizes the pairwise margin for the most risky pair,
                while achieving the least hinge loss
                for general training dataset \(\mathcal{D}\)
                (see Remark~\ref{rem:soft-marginSVM}(b)).
    }
  \end{enumerate}
\end{remark}
\subsection{Fixed-point characterization of \(\mathfrak{S}_\star \)}
To obtain a fixed-point characterization of \(\mathfrak{S}_\star\) in (\ref{Proposition1})
    with computable \(T_{\textrm{DRS}}\) operator  in (\ref{TDRS}),
    we use a convenient expression (see  \cite{chierchia2015proximal})
        \begin{align}
            \Phi_{\mathcal{D}}(\mathbf{p}) = \sum_{i=1}^M h_i(A_{\mathcal{D},i} \mathbf{p})
        \end{align}
    in terms of proximable functions
        \begin{align}
         h_i : \mathbb{R}^K \to \mathbb{R}
             : (\xi_1,\dots,\xi_K) \mapsto \max_{j\in\mathcal{I}_K} \{ r_{\mathcal{D},i}^{(j)} + \xi_j\},
        \end{align}
    where
        \begin{align}
          r_{\mathcal{D},i}^{(j)} :=  
         \begin{cases}
               0\quad \textrm{if}\ j = y_i,\\
               1\quad \textrm{otherwise},
           \end{cases}
       \end{align}
    and  \(A_{\mathcal{D},i} : \mathcal{H}_1 := \mathbb{R}^{(N+1)K} \to \mathbb{R}^K\)
        is defined for \(\mathbf{p} := (\mathbf{w}_j,b_j)_{j\in\mathcal{I}_K}\) as
            \begin{align}      
                A_{\mathcal{D},i} ((\mathbf{w}_j,b_j)_{j\in\mathcal{I}_K} ) :=          
            \left(
            (\mathbf{w}_{y_i}-\mathbf{w}_j)^{\top} \mathbf{x}_i + (b_{y_i}-b_j)
        \right)_{j \in \mathcal{I}_K}. \notag
            \end{align}
To design \(T_{\textrm{DRS}}\) with \(\prox_{h_i}\ (i=1,2,\dots,M) \) as its computable building blocks,
        we  translate a minimization of 
                \(\Phi_\mathcal{D}\) in (\ref{Proposition1})
            into a minimization of sum of 
                \(f \in \Gamma_{0}(\mathcal{X}) \) 
            and
                \(g \in \Gamma_{0}(\mathcal{X})\),
            where
                 \(\mathcal{X} :=\mathcal{H}_1 \times \mathcal{H}_2  := \mathcal{H}_1 \times \mathbb{R}^{KM} \),
                 \( f : \mathcal{X} \to (-\infty,\infty]
                 : (\mathbf{p},\mathbf{u}) :=(\mathbf{p},(\mathbf{u}_1,\dots,\mathbf{u}_M)) \mapsto  \sum_{i=1}^M h_i(\mathbf{u}_i)\)
            and 
                \(g : \mathcal{X} \to (-\infty,\infty]  
                : (\mathbf{p},\mathbf{u}) \mapsto 
                \iota_{\mathcal{N}(\check A_{\mathcal{D}})}(\mathbf{p},\mathbf{u}) \)
            with the null space of
            $\check A_{\mathcal{D}} : \mathcal{X} \to \mathbb{R}^{KM} 
            : (\mathbf{p},\mathbf{u}) \mapsto (A_{\mathcal{D},i}\mathbf{p})_{i=1}^M - \mathbf{u}$.
Since proximity operators of \( f \) and \(g\) 
    are available\footnote{
        \(P_{\mathcal{N}(\check A_{\mathcal{D}})}\) 
            is  computable 
            as a linear operator.
    } as
    \begin{align}
        \left\{
        \begin{array}{ll}
        \prox_{f} (\mathbf{p},\mathbf{u}) 
        = (\mathbf{p},(\prox_{h_1}(\mathbf{u}_1),\dots,\prox_{h_M}(\mathbf{u}_M))); \notag\\
        \prox_{g} (\mathbf{p},\mathbf{u}) = P_{\mathcal{N}(\check A_{\mathcal{D}})}(\mathbf{p},\mathbf{u}) \notag,
        \end{array}
        \right.
    \end{align} 
we obtain, from (\ref{TDRS}) and (\ref{eq:TDRS}),
    \begin{align}
        \hspace{-2mm}
        \mathfrak{S}_{\star}
        &=\mathcal{Q}\left[\argmin (f + g) (\mathcal{X})\right] \notag \\
        &=\mathcal{Q}  P_{\mathcal{N}(\check A_{\mathcal{D}})} (\fix(T_{\textrm{DRS}}))
        \label{TDRS24}
    \end{align}
    with computable \(T_{\textrm{DRS}}\)
    and
        \(\mathcal{Q} :\mathcal{X} \to \mathbb{R}^{(N+1)K} 
                        : (\mathbf{p},\mathbf{u}) \mapsto \mathbf{p}\).
\subsection{How can we achieve rHC-mSVM ?}
Next theorem presents
    a translation of Definition~\ref{def:rHC-mSVM} into a smooth convex optimization over a fixed-point set.
By applying the HSDM (\ref{HSDM}) to this translated problem, 
    we propose an iterative algorithm to solve Problem (\ref{Proposition1}).
\begin{theorem}
Let  \(\mathcal{H}_3 := \mathbb{R}^{\frac{K(K-1)N}{2}} \)
    and
        \( \mathcal{H} := \mathcal{H}_1  \times \mathcal{H}_2 \times \mathcal{H}_3 \times \mathbb{R}\)
    be the Hilbert space where the standard inner product and norm are defined.
Choose \(\Upsilon , \varrho_1 \in \mathbb{R}_{++}\) 
    large enough to satisfy
        \(\|\mathbf{x}_i\| \leq \Upsilon \ (i=1,2,\dots,M)\)
    and
        \(\argmin\Psi(\mathfrak{S}_\star) \cap \bar B_{\mathcal{H}_1}(0,\varrho_1) \neq \varnothing \)
    in Definition~\ref{def:rHC-mSVM}.
Define a  bounded closed convex set 
         \(\mathfrak{B} := [ \times_{i=1}^3 \bar B_{\mathcal{H}_i}(0,\varrho_i) ] \times [-\varrho_4,\varrho_4] \subset \mathcal{H} \),
    where
        \(\varrho_2 :=2\varrho_1\sqrt{KM(\Upsilon^2 + 1)} \),
        \(\varrho_3 := \varrho_1 K(K-1) \)
    and
        \(\varrho_4 := 2\varrho_1\).
    \\[1mm]
{\rm (a)}     \( \mathbf{p}^\diamond \in \argmin \Psi(\mathfrak{S}_\star)  \) (in Def.~\ref{def:rHC-mSVM}),
        where
            \( (\mathbf{p}^\diamond,\mathbf{u}^\diamond,\pmb{\omega}^\diamond,t^\diamond) \in \mathcal{H}\) 
        is a solution of 
            \begin{align}
                    \left.
                    \begin{array}{ll}
                    \minimize 
                    \Theta(\mathbf{p},\mathbf{u},\pmb{\omega},t) := t\  
                    {\rm subject\ to}\ 
            \\\\
            (\mathbf{p},\mathbf{u},\pmb{\omega},t) 
                \in 
            [ \fix(T_{\textrm{DRS}}) \times \mathcal{H}_3 \times \mathbb{R} ]
            \cap \mathfrak{B}
            \\
            \hspace{2cm}
            \cap
            {\displaystyle \bigcap_{r,s \in \mathcal{I}_K\ (r<s)}  } 
            \left( E_{r,s} \cap L_{r,s} \right),
            \end{array}
            \right \}
            \label{MSVM4}
            \end{align}
        where \(T_{\textrm{DRS}}\) is defined in (\ref{TDRS24}),
            \( E_{r,s} := \{ (\mathbf{p},\mathbf{u},\pmb{\omega},t) \in \mathcal{H} \mid  \| \pmb{\omega}_{rs}\| \leq t\} \),
        and
            \( L_{r,s} := \{(\mathbf{p},\mathbf{u},\pmb{\omega},t) \in \mathcal{H} \ |\  \pmb{\omega}_{rs} = \Omega_{r,s}P_{\mathcal{N}(\check A_\mathcal{D})} (\mathbf{p},\mathbf{u}) \} \)
        for
            \(\pmb{\omega} = (\pmb{\omega}_{r's'})_{r'<s'} \in \mathcal{H}_3 \)
        with
        $\Omega_{(r,s)}  : \mathcal{H}_1 \times \mathcal{H}_2 \to \mathbb{R}^N : (\mathbf{p},\mathbf{u}) := ((\mathbf{w}_j,b_j)_{j \in \mathcal{I}_K},\mathbf{u})
    \mapsto \textbf{w}_r - \textbf{w}_s $.
    \\[1mm]    
{\rm (b)} Define for \(\alpha \in (0,1) \) 
        \begin{align}
            T :=
            P_{\mathfrak{B}} \circ \widetilde T_{\textrm{DRS}_\alpha}\circ &P_{E_{K-1,K}}\circ \cdots \notag
                \circ P_{E_{1,2}}\\ &\circ P_{L_{K-1,K}} \circ \cdots \circ P_{L_{1,2}} \label{T},
        \end{align}    
    where
            \( \widetilde T_{\textrm{DRS}_\alpha}  : \mathcal{H} \to \mathcal{H}
                                    : (\mathbf{p},\mathbf{u},\pmb{\omega},t) 
                                    \mapsto ((1-\alpha) (\mathbf{p},\mathbf{u})+ \alpha T_{\textrm{DRS}} (\mathbf{p},\mathbf{u}),\pmb{\omega},t)
            \).
Then \(T\) is a computable averaged nonexpansive operator  
    with        
        \begin{align}
            \fix(T)
            =
            [ \fix(T_{\textrm{DRS}}) \times \mathcal{H}_3 \times \mathbb{R} ] 
            \cap
            \mathfrak{B}\notag \\
            \cap
           {\displaystyle \bigcap_{r,s \in \mathcal{I}_K\ (r<s)}  } 
        \left( E_{r,s} \cap L_{r,s} \right)
        \label{FixT}
        \end{align}
        (Note: See \cite{bauschke1996projection1} for an explicit expression of \(P_{E_{r,s}}\) in (\ref{T})).
    \\[1mm]
{\rm (c)} By using  $\Theta$ in (\ref{MSVM4}) and $T$ in (\ref{T}), 
    the sequence \(\mathbf{z}_{n} := (\mathbf{p}_n,\mathbf{u}_n,\pmb{\omega}_n,t_n) \in \mathcal{H} \ (n=0,1,2\dots)\),
        generated with any \(\mathbf{z}_0 \in \mathcal{H}\) by the HSDM (\ref{HSDM}), satisfies
            \begin{align}
                \left \{
                \begin{array}{ll}
                    {\displaystyle \lim_{n \to \infty}
                \min_{\mathbf{z}^\diamond \in \mathfrak{S}_{\mathcal{H}}^\diamond} 
                \|\mathbf{z}_n - \mathbf{z}^\diamond \|_{\mathcal{H}} = 0};
        \\
        {\displaystyle   
            \lim_{n \to \infty}
                \min_{\mathbf{p}^{\star\star} \in \argmin\Psi(\mathfrak{S}_\star) \cap \bar B_{\mathcal{H}_1}(0,\varrho_1)} 
                \|\mathbf{p}_n - \mathbf{p}^{\star\star}\|_{\mathcal{H}_1} = 0,
                    }
                \end{array}
                \right.\notag
            \end{align}
    where
        \(\mathfrak{S}_{\mathcal{H}}^\diamond 
        := \{\mathbf{z} \in \mathcal{H} \mid \mathbf{z} \textrm{\ is a solution of\ } (\ref{MSVM4}) \} \).
      \end{theorem}
\section{NUMERICAL EXPERIMENTS}
In our experiment,
    we applied the naive convex relaxation in (\ref{pro:Crammer_Singer}) 
        and the proposed hierarchical formulation in (\ref{Proposition1})
            to the Iris dataset used in \cite{fisher1936use}.
This data \(\mathcal{D}\) has 150 labeled sample points, 
    which are divided into three classes 
        \(\mathcal{I}_3 := \{1(\textrm{setosa}), 2(\textrm{versicolor}), 3(\textrm{virginica})  \}\),
    i.e.,
        \(\mathcal{D}= \{((x_i^{SL},x_i^{SW},x_i^{PL},x_i^{PW}),y_i) \in \mathbb{R}^4\times \mathcal{I}_3 \mid i=1,\dots,150 \} \),
    where \(x_i^{SL}\), \(x_i^{SW}\), \(x_i^{PL}\) and \(x_i^{SL}\)
        stand respectively for sepal length, sepal width, petal length, and petal width.
 We chose \(\mathfrak{I}_{\textrm{sep}},\mathfrak{I}_{\textrm{nsep}} \subset \{1,\dots,150\}\ ( |\mathfrak{I}_{\textrm{sep}} | =  |\mathfrak{I}_{\textrm{nsep}} | = 75 )\),
    to make 
        \(\mathcal{D}_{\textrm{sep}}= \{((x_i^{SW},x_i^{PW}),y_i) \in \mathbb{R}^2\times \mathcal{I}_3 \mid i\in\mathfrak{I}_{\textrm{sep}} \} \)
    linearly separable and
        \(\mathcal{D}_{\textrm{nsep}}= \{((x_i^{SL},x_i^{SW}),y_i) \in \mathbb{R}^2\times \mathcal{I}_3 \mid i\in\mathfrak{I}_{\textrm{nsep}} \} \)
    linearly non-separable.

\begin{figure}[t]

    \begin{minipage}[b]{1\linewidth}
        \centering
        \centerline{\includegraphics[width=10cm]{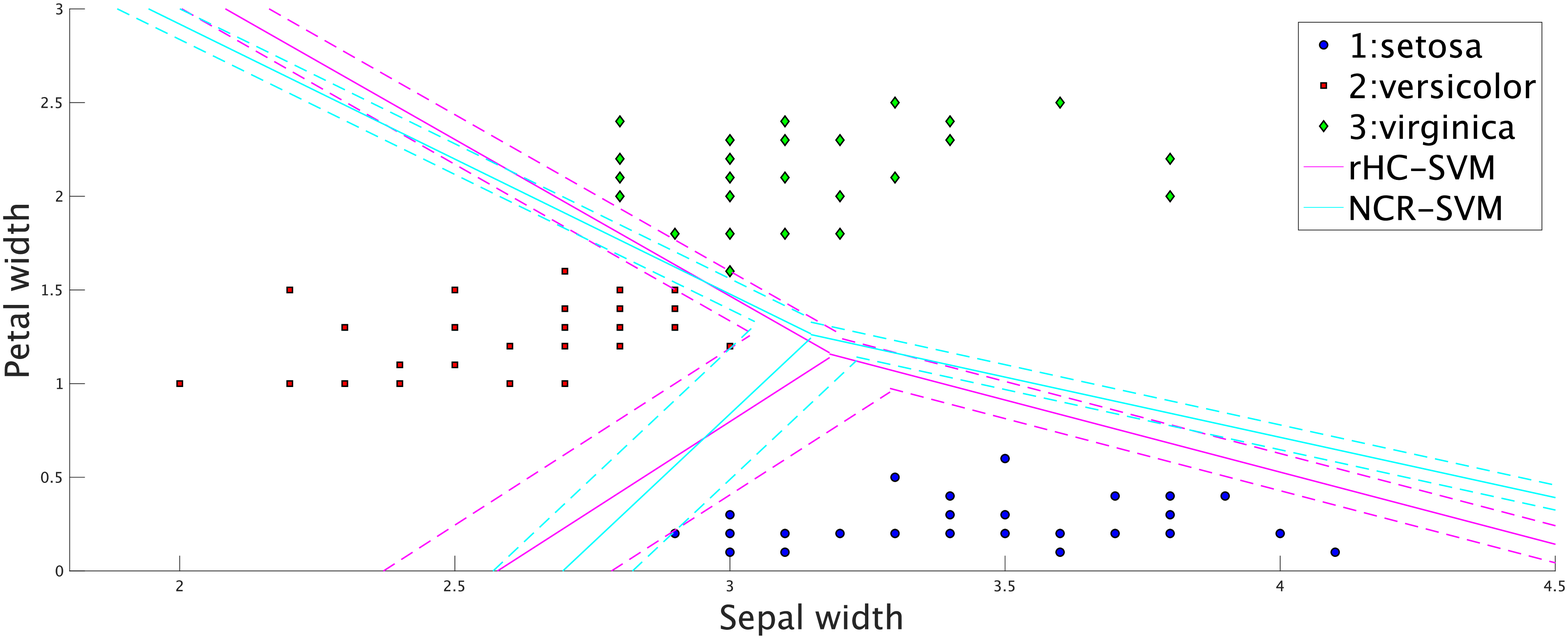}}
        \caption{Experiment for separable  dataset \(\mathcal{D}_{\textrm{sep}}\)}\medskip
        \label{fig:ressep}
      \end{minipage}
      \begin{minipage}[b]{1\linewidth}
        \centering
        \centerline{\includegraphics[width=10cm]{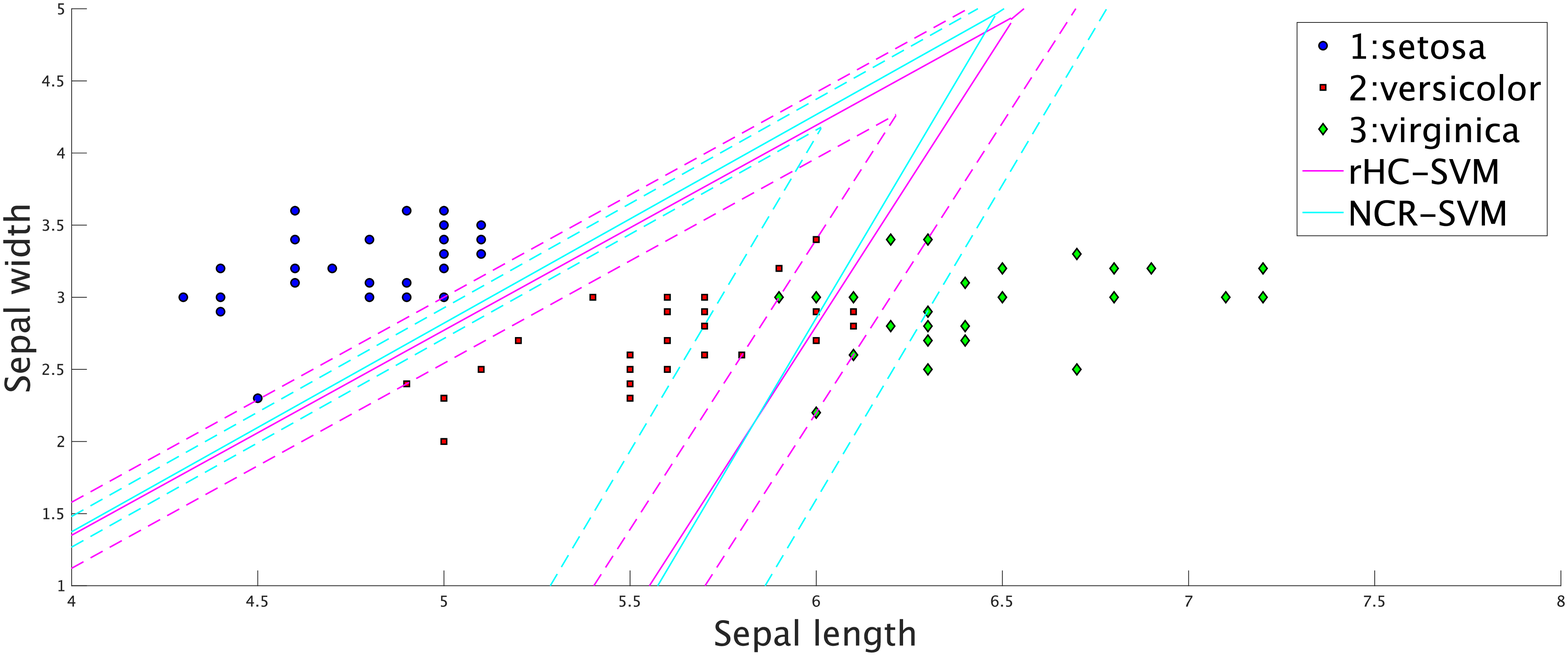}}
        \caption{Experiment for non-separable  dataset \(\mathcal{D}_{\textrm{nsep}}\)}\medskip
        \label{fig:resnsep}        
      \end{minipage}
    %
\end{figure}
To see the performance of the linear classifier as the solution of (\ref{pro:Crammer_Singer}),
    we applied a Forward-Backward based Primal-Dual method \cite{vu2013splitting,condat2013primal,chierchia2015proximal}
        with \( C = 2^{10}\) to \(\mathcal{D}_{\textrm{sep}}\) and \(\mathcal{D}_{\textrm{nsep}}\).
    
To see the performance of the linear classifier as the solution of (\ref{Proposition1}),
    we applied the proposed algorithmic solution in Theorem 1(c)
        to \(\mathcal{D}_{\textrm{sep}}\) and \(\mathcal{D}_{\textrm{nsep}}\).

For each linear classifier
    \(\mathbf{p} = \mathbf{p}_{\textrm{NCR}}\) by (\ref{pro:Crammer_Singer}) and     
    \(\mathbf{p} = \mathbf{p}_{\textrm{rHC}}\) by (\ref{Proposition1}),
    the three decision boundaries 
        \(\Pi_{(\mathbf{p},\textrm{r,s})}^{0} (r<s)\)
    are drawn with solid lines in  Fig.~\ref{fig:ressep} for 
        \(\mathcal{D}_{\textrm{sep}}\)
    and in Fig.~\ref{fig:resnsep} for 
        \(\mathcal{D}_{\textrm{nsep}}\).
The boundaries of \(\Pi_{(\mathbf{p},\textrm{r,s})}^{[1,\infty)} (r \neq s)\)     
    are drawn with dashed lines in Fig.~\ref{fig:ressep} for \(\mathcal{D}_{\textrm{sep}}\) and Fig.~\ref{fig:resnsep} for \(\mathcal{D}_{\textrm{nsep}}\).    

In Fig.~\ref{fig:ressep},
    we observe that both linear classifiers achieve 
    \(\mathscr{E}_{\mathbf{p}_{\textrm{NCR}}} (\mathcal{D}_{\textrm{sep}})
        = \mathscr{E}_{\mathbf{p}_{\textrm{rHC}}}(\mathcal{D}_{\textrm{sep}}) = 0\)
    with smallest pairwise margins
    \(\textrm{dist} \left(\Pi_{(\mathbf{p}_{\textrm{NCR}},2,3)}^{[1,\infty)} ,\Pi_{(\mathbf{p}_{\textrm{NCR}},2,3)}^{0}  \right) = 0.0467\)
    and\\
    \(\textrm{dist} \left(\Pi_{(\mathbf{p}_{\textrm{rHC}},2,3)}^{[1,\infty)} ,\Pi_{(\mathbf{p}_{\textrm{rHC}},2,3)}^{0}  \right) = 0.0681\)
(Note:
    The solutions \(\mathbf{p}_{\textrm{NCR}}\) of 
    (\ref{pro:Crammer_Singer}) with \(C=2^8 \ \textrm{and} \ 2^9\) result respectively in 
    \(\mathscr{E}_{\mathbf{p}_{\textrm{NCR}}} (\mathcal{D}_{\textrm{sep}})=3 \ \textrm{and}\  1 \)).

In Fig.2, 
    we observe that both linear classifiers achieve 
        \(\mathscr{E}_{\mathbf{p}_{\textrm{NCR}}} (\mathcal{D}_{\textrm{nsep}}) = 16\)
    and 
        \( \mathscr{E}_{\mathbf{p}_{\textrm{rHC}}}(\mathcal{D}_{\textrm{nsep}}) = 7\)
    with smallest pairwise margins
        \(\textrm{dist} \left(\Pi_{(\mathbf{p}_{\textrm{NCR}},1,3)}^{[1,\infty)} ,\Pi_{(\mathbf{p}_{\textrm{NCR}},1,3)}^{0}  \right) = 0.0503\)
    and\\
        \(\textrm{dist} \left(\Pi_{(\mathbf{p}_{\textrm{rHC}},1,3)}^{[1,\infty)} ,\Pi_{(\mathbf{p}_{\textrm{rHC}},1,3)}^{0}  \right) = 0.0702\).

These experiments demonstrate that
    the desired features in Remark~\ref{rem:rHC-mSVM} are achieved by the proposed algorithm in Theorem 1(c).
We also observed the difficulty in parameter tuning for (\ref{pro:Crammer_Singer}).

\section{CONCLUSION}
We proposed the Robust Hierarchical Convex multiclass SVM and 
an iterative algorithm to realize this classifier.
Numerical experiments demonstrate that 
the proposed multiclass classifier achieves
larger smallest pairwise margin, even with smaller number of misclassified training samples, than a multiclass classifier [Crammer-Singer 2001]
which has been known as a most natural multiclass extension of the soft-margin SVM.

%

\balance
\bibliographystyle{IEEEbib}
\bibliography{EUSIPCO_NYY-2020j}

\end{document}